# UNIFORM ASYMPTOTICS FOR ROBUST LOCATION ESTIMATES WHEN THE SCALE IS UNKNOWN[1]

By Matias Salibian-Barrera and Ruben H. Zamar

*Carleton University and University of British Columbia*

Most asymptotic results for robust estimates rely on regularity conditions that are difficult to verify in practice. Moreover, these results apply to fixed distribution functions. In the robustness context the distribution of the data remains largely unspecified and hence results that hold uniformly over a set of possible distribution functions are of theoretical and practical interest. Also, it is desirable to be able to determine the size of the set of distribution functions where the uniform properties hold. In this paper we study the problem of obtaining verifiable regularity conditions that suffice to yield uniform consistency and uniform asymptotic normality for location robust estimates when the scale of the errors is unknown. We study $M$-location estimates calculated with an $S$-scale and we obtain uniform asymptotic results over contamination neighborhoods. Moreover, we show how to calculate the maximum size of the contamination neighborhoods where these uniform results hold. There is a trade-off between the size of these neighborhoods and the breakdown point of the scale estimate.

**1. Introduction.** Many robust location point estimates have been proposed in the last 35 years. Unfortunately, robust inference has not received the same amount of attention in the literature. Since the finite sample distributions of robust estimates are unknown, robust inference typically relies on the asymptotic distributions of these estimates.

According to the robustness model one does not know the actual distribution of the data. Therefore it is highly desirable to have asymptotic results that hold uniformly over some set of plausible distributions. Hampel (1971) considered this problem and broke new ground showing that under certain regularity conditions $M$-location estimates have uniform asymptotic

Received July 2002; revised April 2003.
[1]Supported by NSERC through individual grants.
*AMS 2000 subject classifications.* 62F35, 62F12, 62E20.
*Key words and phrases.* Robustness, robust inference, robust location and scale models, $M$-estimates.

properties on Prokhorov neighborhoods. These are the first results in the robustness literature that deal with uniform asymptotic properties of robust estimates. Moreover, Hampel's results provide a valuable set of tools to evaluate and compare robust estimates based on their asymptotic behavior. Unfortunately his results only guarantee the existence of a neighborhood with unknown size where the uniform behavior holds. Determining the size of these neighborhoods seems to be a very difficult problem [see also Davies (1998)].

Huber [(1967); (1981), page 51] shows that when the scale of the errors is known, the $M$-location estimates are asymptotically normal and the approximation is uniform on the set of symmetric distributions that have all their mass concentrated on the points where the estimating equation is differentiable. Simultaneous estimation of location and scale with Huber's Proposal 2 was studied by Clarke (1980, 1986). In the first reference the author considers the problem of uniform convergence for these estimates and in the second proves that Huber's Proposal 2 estimates with nonsmooth estimating equations fall in the framework of Hampel (1971). More recently, Davies (1998) constructed $M$-location estimates with simultaneous scale estimates (Huber's Proposal 2) that are locally asymptotically normal. Davies's results are "locally uniform"; that is, for each distribution function there exists a neighborhood of distributions where the convergence holds uniformly. Unfortunately, the size of these neighborhoods is unknown. Finally, Clarke (2000) shows that certain $M$-location estimates [including the simultaneous location and scale estimation proposed in Heathcote and Silvapulle (1981)] are continuous over full Prokhorov neighborhoods of the parametric model. It follows that these estimates have uniform asymptotic behavior over Prokhorov neighborhoods. Unfortunately, as in Hampel (1971) and Davies (1998), the size of these neighborhoods is unknown.

Our results apply to location $M$-estimates calculated using an $S$-scale [see Rousseeuw and Yohai (1984)]. These estimates are scale-equivariant [Rousseeuw and Leroy (1987), pages 158 and 159] and have simultaneous high breakdown point and high efficiency at the central model. Moreover, under verifiable regularity conditions we obtain uniform asymptotic results (consistency and asymptotic distribution) that hold over a contamination neighborhood of *known size*. The size of these sets depends on the breakdown point of the $S$-scale estimates (the higher the breakdown point, the smaller the set of distribution functions where uniformity holds; see Table 1, in Section 4). The regularity conditions we need in our results depend on two separate aspects of the inference procedure: the assumed "true" parametric model for the "good" data points, and the estimating equations used to calculate the robust estimate.

The rest of the paper is organized as follows. Section 2 contains the definitions of the estimates we consider. Section 3 shows that under mild regularity conditions these estimates are uniformly consistent on contamination neighborhoods. Section 4 gives additional assumptions under which the above estimates are uniformly asymptotically normal. Section 5 contains some concluding remarks and, finally, the Appendix contains sketches of the proofs of our main results. Details can be found in Salibian-Barrera and Zamar (2002).

**2. MM-location estimates.** Consider the following location–scale model: let $x_1, \ldots, x_n$ be $n$ observations on the real line satisfying

$$x_i = \mu + \sigma \varepsilon_i, \qquad i = 1, \ldots, n, \tag{1}$$

where $\varepsilon_i$, $i = 1, \ldots, n$, are independent and identically distributed (i.i.d.) observations with variance equal to 1. The interest is in estimating $\mu$, and the scale $\sigma$ is considered a nuisance parameter.

We will consider scale-equivariant $M$-location estimates $\hat{\mu}_n$ defined as the solution of an estimating equation of the form

$$\frac{1}{n} \sum_{i=1}^{n} \psi((x_i - \hat{\mu}_n)/\hat{\sigma}_n) = 0, \tag{2}$$

where $\hat{\sigma}_n$ is an $S$-scale estimate of the residuals [Rousseeuw and Yohai (1984)] and $\psi: \mathbb{R} \to \mathbb{R}$ is a nondecreasing, odd and continuously differentiable real function. An example of such a function is given by

$$\psi_c(u) = \text{sign}(u) \begin{cases} |u/c|, & \text{if } |u| \leq 0.8c, \\ p_4(|u|/c), & \text{if } 0.8c < |u| \leq c, \\ p_4(1), & \text{if } |u| > c, \end{cases} \tag{3}$$

where $c > 0$ is a user-chosen tuning constant, and $p_4(u) = 38.4 - 175u + 300u^2 - 225u^3 + 62.5u^4$ [see Fraiman, Yohai and Zamar (2001) and also Bednarski and Zontek (1996), for other choices of smooth functions $\psi$]. Following Yohai (1987) we will call these $M$-location estimates obtained with an $S$-scale *MM-location* estimates.

The $S$-scale estimate $\hat{\sigma}_n$ we use in (2) is defined as follows. Let $\rho: \mathbb{R} \to \mathbb{R}_+$ be a bounded, continuous and even function satisfying $\rho(0) = 0$. The $S$-scale $\hat{\sigma}_n$ is defined by

$$\hat{\sigma}_n = \inf_{t \in \mathbb{R}} s_n(t), \tag{4}$$

where, for each $t \in \mathbb{R}$, $s_n(t)$ is the solution of

$$\frac{1}{n} \sum_{i=1}^{n} \rho((x_i - t)/s_n(t)) = b, \tag{5}$$

and $b = E[\rho(u)]$. Naturally associated with this family are the *S-location* estimates $\tilde{\mu}_n$ given by

$$\tilde{\mu}_n = \arg\inf_{t \in \mathbb{R}} s_n(t). \tag{6}$$

Beaton and Tukey (1974) proposed a family of functions $\rho_d$ given by

$$\rho_d(u) = \begin{cases} 3(u/d)^2 - 3(u/d)^4 + (u/d)^6, & \text{if } |u| \leq d, \\ 1, & \text{if } |u| > d, \end{cases} \tag{7}$$

where the tuning constant $d$ is positive. The above family of functions $\rho_d$ satisfies all the regularity conditions we need to obtain uniform asymptotic properties, and at the same time it yields scale estimates $\hat{\sigma}_n$ with good robustness properties.

REMARK 1 ($\psi \neq \rho'$). Note that the estimating function $\psi$ in (2) need not be equal to $\rho'$ in (5). Moreover, we will recommend using $\psi = \psi_c$ in (3) and $\rho = \rho_d$ in (7).

REMARK 2 (High efficiency and high breakdown point). The robust location estimates $\hat{\mu}_n$ defined by (2) with $\hat{\sigma}_n$ as in (4) are scale-equivariant and can have simultaneously high breakdown and high efficiency at the central model. For example, the choice $d = 1.548$ for $\rho_d$ in (7), $b = 0.5$ in (5) and $c = 1.525$ for $\psi_c$ in (3) yields a location estimate $\hat{\mu}_n$ with 50% breakdown point and 95% efficiency when the errors have a normal distribution.

The asymptotic properties (consistency and asymptotic normality) of $M$-location estimates given by (2) are well known when the distribution of the errors is symmetric [Huber (1964, 1967, 1981), Boos and Serfling (1980) and Clarke (1983, 1984)]. The next two sections establish these properties uniformly over a set of distributions.

**3. Uniform consistency.** The goal of this section is to determine verifiable conditions under which the *scale-equivariant* $M$-location estimates $\hat{\mu}_n$ given by (2) are uniformly consistent on the contamination "neighborhood"

$$\mathcal{H}_\varepsilon(F_0) = \{F \in \mathcal{D} : F(x) = (1 - \varepsilon) F_0((x - \mu)/\sigma) + \varepsilon H(x)\}, \tag{8}$$

where $\mathcal{D}$ denotes the set of all distribution functions, $F_0$ is a fixed symmetric distribution, $\mu_0$ and $\sigma_0$ are the unknown location and scale parameters, $\varepsilon \in (0, 1/2)$ and $H$ is an arbitrary distribution function. Since in what follows the central distribution $F_0$ is fixed, we write $\mathcal{H}_\varepsilon$ to denote the set (8) above.

Under certain regularity conditions (see references above) the $M$-location estimates $\hat{\mu}_n$ and the $S$-estimates $\hat{\sigma}_n$ and $\tilde{\mu}_n$ are consistent to the functionals $\boldsymbol{\mu}(F)$, $\boldsymbol{\sigma}(F)$ and $\tilde{\boldsymbol{\mu}}(F)$ defined by the following equations. For each $t \in \mathbb{R}$, let $\sigma(F,t)$ satisfy

$$E_F[\rho((X-t)/\sigma(F,t))] = b. \tag{9}$$

The asymptotic value of $\hat{\sigma}_n$ is given by

$$\boldsymbol{\sigma}(F) = \inf_{t \in \mathbb{R}} \sigma(F,t). \tag{10}$$

Similarly, for the $S$-location estimate $\tilde{\mu}_n$ we have

$$\tilde{\boldsymbol{\mu}}(F) = \arg\inf_{t \in \mathbb{R}} \sigma(F,t). \tag{11}$$

Finally, for the $M$-location estimate $\hat{\mu}_n$ the corresponding equation is

$$E_F[\psi((X - \boldsymbol{\mu}(F))/\boldsymbol{\sigma}(F))] = 0. \tag{12}$$

DEFINITION 1 (Uniform consistency). We say that the sequence of estimates $\hat{\tau}_n$ is uniformly consistent to the functional $\tau(F)$ over the contamination neighborhood $\mathcal{H}_\varepsilon$ if, for all $\delta > 0$,

$$\lim_{m \to \infty} \sup_{F \in \mathcal{H}_\varepsilon} P_F\left[\sup_{n \geq m} |\hat{\tau}_n - \tau(F)| > \delta\right] = 0,$$

where $\tau(F)$ is the a.s. limit of $\hat{\tau}_n$ for an i.i.d. sequence of observations with distribution function $F$. We will denote this type of convergence by $\hat{\tau}_n \xrightarrow{\varepsilon} \tau$.

Our main result in this section states that if the scale estimate $\hat{\sigma}_n$ in (2) satisfies $\hat{\sigma}_n \xrightarrow{\varepsilon} \sigma$ and if $\psi$ is odd, nondecreasing, bounded and continuously differentiable, then $\hat{\mu}_n \xrightarrow{\varepsilon} \mu$.

THEOREM 1 (Uniform consistency of the $M$-location estimate with general scale). Let $x_1, \ldots, x_n$ be i.i.d. observations following the location model (1). Let $\psi$ satisfy the following:

(P.1) $|\psi(u)| \leq 1$ for all $u \in \mathbb{R}$, and $\psi(-u) = -\psi(u)$ for $u \geq 0$;
(P.2) $\psi$ is nondecreasing and $\lim_{u \to \infty} \psi(u) > 0$;
(P.3) $\psi$ is continuously differentiable.

Suppose that $\hat{\sigma}_n$ in (2) has asymptotic breakdown point $\varepsilon^*$. Let $0 \leq \varepsilon < \varepsilon^*$ be such that $\hat{\sigma}_n \xrightarrow{\varepsilon} \sigma$. Then if $\hat{\mu}_n$ satisfies (2), we have $\hat{\mu}_n \xrightarrow{\varepsilon} \mu$.

A sketch of the proof of Theorem 1 is given in the Appendix. A detailed proof can be found in Salibian-Barrera and Zamar (2002).

REMARK 3 (Uniform consistency of $S$-scale estimates). When $\hat{\sigma}_n$ is an $S$-scale estimate, Martin and Zamar (1993) showed that if $F_0$ (the central distribution function in $\mathcal{H}_\varepsilon$) has an even and unimodal density, and if the function $\rho$ is even, bounded, continuous and nondecreasing in $[0, \infty)$, then $\hat{\sigma}_n$ has asymptotic breakdown point $1/2$. They also showed that if in addition $F_0$ has a positive density on the real line, then for all $0 < \varepsilon < 1/2$ we have

$$\hat{\sigma}_n \xrightarrow{\varepsilon} \sigma. \tag{13}$$

Theorem 1 and Remark 3 imply that $M$-location estimates $\hat{\mu}_n$ given by (2) with $\psi = \psi_c$ in the family (3) and scale $\hat{\sigma}_n$ given by (4) with $\rho = \rho_d$ in Tukey's family (7) have high breakdown point and high efficiency and are uniformly consistent over $\mathcal{H}_\varepsilon$ for all $0 < \varepsilon < 1/2$.

**4. Uniform asymptotic distribution.** In this section we show that under certain regularity conditions the $MM$-location estimates $\hat{\mu}_n$ converge weakly to a normal distribution uniformly over the contamination neighborhood $\mathcal{H}_\varepsilon$. These results are constructive and allow us to determine the size of the neighborhood $\mathcal{H}_\varepsilon$ where uniform asymptotic normality holds. The required regularity conditions will be mainly imposed on our estimating equations (2) and (4) and we will show that $\psi = \psi_c$ in (3) and $\rho = \rho_d$ in (7) satisfy these conditions. Hence, our results show that the scale-equivariant $MM$-location estimates have simultaneously high breakdown point, high efficiency at the central model and are uniformly asymptotically normal on a contamination neighborhood of known size (see Remark 2).

Asymptotic results for asymmetric distributions are not easy to obtain. There are some results in the robustness literature dealing with this problem [Carroll (1978, 1979), Carroll and Welsh (1988) and Rocke and Downs (1981)]. They show that when $F$ is asymmetric the asymptotic distribution of the location estimate depends on that of the scale and that the asymptotic variance calculated with the assumption of symmetry is not correct. Salibian-Barrera (2000) showed that in general the asymptotic distribution of location $M$-estimates for arbitrary distribution functions when the scale is estimated with an $S$-scale depends on the behavior of the $S$-scale and the corresponding $S$-location estimate as well. Hence, to obtain uniform asymptotics for these $MM$-location estimates we need uniform consistency of the *S-scale* and *S-location estimates*.

$S$-scale estimates are uniformly consistent under relatively weak regularity conditions [see Martin and Zamar (1993) and our Remark 3].

Uniform consistency of *S-location* estimates requires more assumptions. For a given $0 \leq \varepsilon < 1/2$ and an estimating function $\rho$ in (5) let $s^+$ and $s^-$ satisfy

$$0 < s^- \leq \inf_{F \in \mathcal{H}_\varepsilon} \boldsymbol{\sigma}(F) < \sup_{F \in \mathcal{H}_\varepsilon} \boldsymbol{\sigma}(F) \leq s^+ < \infty. \tag{14}$$

Note that, from Lemma 1 in Martin and Zamar (1993), for all $0 \leq \varepsilon < \min(b, 1-b)$ we have $0 < \inf_{F \in \mathcal{H}_\varepsilon} \boldsymbol{\sigma}(F) \leq \sup_{F \in \mathcal{H}_\varepsilon} \boldsymbol{\sigma}(F) < \infty$. To simplify the notation we will omit the dependence of $s^+$ and $s^-$ on $\varepsilon$. Assume that there exists $t^* \in \mathbb{R}$ such that

$$(15) \qquad \inf_{s^- \leq s \leq s^+} \left[ E_{F_0} \rho\left(\frac{X-t}{s}\right) - E_{F_0} \rho\left(\frac{X}{s}\right) \right] > \frac{\varepsilon}{1-\varepsilon} \qquad \forall |t| \geq t^*$$

and

$$(16) \qquad \inf_{-t^* \leq t \leq t^*, s^- \leq s \leq s^+} E_{F_0} \rho''\left(\frac{X-t}{s}\right) > \frac{\varepsilon}{1-\varepsilon} \sup_x [\rho''(x)]^-,$$

where $s^+$ and $s^-$ are given in (14).

Condition (16) can be slightly relaxed [see Salibian-Barrera and Zamar (2002), Lemma 7]. Assumptions (15) and (16) do not depend on the unknown distribution of the data $F$ (only on $F_0$, the central distribution of the neighborhood $\mathcal{H}_\varepsilon$) but are tedious to verify and will typically require numerical computations. Note that for a particular $\rho$ these conditions impose an upper bound $\varepsilon = \varepsilon(\rho, F_0)$ on the size of the contamination neighborhood $\mathcal{H}_\varepsilon$. More specifically, fix the central distribution $F_0$ and note that the left-hand side of (15) is a nondecreasing function of $|t|$. At the same time, the right-hand side is an increasing function of $\varepsilon$. Hence, the smallest value $t^* = t^*(\varepsilon)$ that satisfies (15) is a nondecreasing function of $\varepsilon$. Also, the left-hand side of (16) is a nonincreasing function of $t^*$ (and thus of $\varepsilon$), while its right-hand side is increasing in $\varepsilon$. It follows that there is a critical value $\varepsilon(\rho)$ such that both (15) and (16) hold for $\varepsilon \leq \varepsilon(\rho)$, but fail to hold for $\varepsilon > \varepsilon(\rho)$. When $\rho = \rho_d$ belongs to Tukey's family (7) and the center of the contamination neighborhood is the standard normal distribution $F_0 = \Phi$ we used numerical methods to find $\varepsilon(\rho_d)$ for different choices of the tuning constant $d$ (i.e., for different breakdown points). We found that there is a trade-off between the breakdown point of the scale estimate and the upper bound $\varepsilon(\rho_d)$: the larger the breakdown point, the smaller the upper bound $\varepsilon(\rho_d)$. Table 1 lists the values of $\varepsilon(\rho_d)$ for contamination neighborhoods of the standard normal distribution and estimating equations that yield estimates with breakdown points between 0.25 and 0.50.

The following theorem states that under these conditions $S$-location estimates are uniformly consistent. This result will be necessary to obtain uniform asymptotic distribution of the $M$-location estimate calculated with an $S$-scale as in (2).

THEOREM 2 (Uniform consistency of the $S$-location estimate). *Suppose that the nonconstant function $\rho$ satisfies the following assumptions:*

(R.1) $\rho(-u) = \rho(u)$, $u \geq 0$, and $\sup_{u \in \mathbb{R}} \rho(u) = 1$;

(R.2) $\rho(u)$ is nondecreasing in $u \geq 0$;
(R.3) $|\rho'(u)| \leq K < \infty$, $\forall\, u \in \mathbb{R}$;
(R.4) there exists $0 < c < \infty$ such that $\rho(u) = 1$ $\forall\, |u| \geq c$.

Let $b \in (0,1)$, $\tilde{\mu}_n$ as in (6) and $\tilde{\boldsymbol{\mu}}(F)$ as in (11). Let $s^+$ and $s^-$ be as in (14) and suppose that $0 < \varepsilon$ is such that (15) and (16) hold. Then

$$(17) \qquad \lim_{m \to \infty} \sup_{F \in \mathcal{H}_\varepsilon} P_F\bigg( \sup_{n \geq m} |\tilde{\mu}_n - \boldsymbol{\mu}(F)| > \delta \bigg) = 0.$$

A sketch of the proof of Theorem 2 is given in the Appendix. A detailed proof can be found in Salibian-Barrera and Zamar (2002).

We can now state our main result: when the $M$-location, $S$-scale and $S$-location estimates are uniformly consistent, the $M$-location estimate has a uniformly asymptotically normal distribution.

THEOREM 3 (Uniform asymptotic distribution of $MM$-location estimates). Let $\hat{\mu}_n$ satisfy (2) with a function $\psi$ that satisfies assumptions (P.1) and (P.2) in Theorem 1 and the following:

(P.4) $\psi$ is twice continuously differentiable;
(P.5) there exists $d > 0$ such that $|\psi(u)| = 1$ for all $|u| \geq d$.

Assume that the $S$-scale estimate $\hat{\sigma}_n$ in (2) is given by (4) with a function $\rho$ that satisfies (R.1)–(R.4) in Theorem 2 and the following:

(R.5) $\rho$ is twice continuously differentiable.

Suppose that $\varepsilon$ is such that (15) and (16) hold and that the center $F_0$ of the contamination neighborhood $\mathcal{H}_\varepsilon$ has a positive, even and unimodal density.

TABLE 1
*Maximum size $\varepsilon(d)$ of contamination neighborhoods around the standard normal distribution where uniform consistency of the $S$-location estimate holds for different breakdown points (BP); the column labeled $d$ contains the tuning constant that yields the respective BP*

| BP | $d$ | $\varepsilon(d)$ |
|---|---|---|
| 0.50 | 1.548 | 0.11 |
| 0.45 | 1.756 | 0.14 |
| 0.40 | 1.988 | 0.17 |
| 0.35 | 2.252 | 0.20 |
| 0.30 | 2.561 | 0.24 |
| 0.25 | 2.937 | 0.25 |

*Then*

$$\lim_{n\to\infty} \sup_{F\in\mathcal{H}_\varepsilon} \sup_{x\in\mathbb{R}} \left| P_F\left\{ \sqrt{n}\frac{(\hat{\mu}_n - \mu)}{\sqrt{V}} < x \right\} - \Phi(x) \right| = 0,$$

*where*

$$\begin{aligned}
V &= V(\mu, \sigma, F) \\
&= \sigma(F)^2 H(F)^2 E_F\left\{ \left[ \psi\left(\frac{X - \mu(F)}{\sigma(F)}\right) \right.\right. \\
&\qquad\qquad\qquad\qquad \left.\left. - J(F)\left(\rho\left(\frac{X - \tilde{\mu}(F)}{\sigma(F)}\right) - b\right) \right]^2 \right\},
\end{aligned} \tag{18}$$

$$H(F) = \frac{1}{E_F\{\psi'(t(X - \mu(F))/\sigma(F))\}},$$

*and*

$$J(F) = \frac{E_F\{\psi'((X - \mu(F))/\sigma(F))(X - \mu(F))/\sigma(F)\}}{E_F\{\rho'((X - \tilde{\mu}(F))/\sigma(F))(X - \tilde{\mu}(F))/\sigma(F)\}}.$$

A sketch of the proof of Theorem 3 is given in the Appendix. A detailed proof can be found in Salibian-Barrera and Zamar (2002).

REMARK 4 (Regularity conditions). The assumptions on $F_0$ (the center of the contamination neighborhood) are needed to show that the $S$-scale estimate $\hat{\sigma}_n$ is uniformly consistent ($\hat{\sigma}_n \xrightarrow{\varepsilon} \sigma$). By Theorem 1 we also have that the $MM$-location estimates are uniformly consistent as well ($\hat{\mu}_n \xrightarrow{\varepsilon} \mu$). The assumptions on the estimating equation $\rho$ of the $S$-scale $\hat{\sigma}_n$ and conditions (15) and (16) are needed to obtain uniform consistency of the $S$-location estimate ($\tilde{\mu}_n \xrightarrow{\varepsilon} \tilde{\mu}$). See Theorem 2.

Using Table 1 we find, for example, that scale-equivariant $MM$-location estimates calculated with $\psi = \psi_{1.525}$ in (3) and an $S$-scale with $\rho = \rho_{1.548}$ in (7) have simultaneously breakdown point $1/2$, have 95% efficiency when the errors are normally distributed, and are uniformly asymptotically normal on a contamination neighborhood of size at least $\varepsilon = 0.11$. If, on the other hand, we use $\rho = \rho_{2.937}$ in (7) we obtain estimates that have the same efficiency, that have lower breakdown point (25%) and that are uniformly asymptotically normal on a contamination neighborhood of size $\varepsilon = 0.25$.

**5. Conclusions.** There are four important properties of robust location estimates: translation and scale-equivariance [Rousseeuw and Leroy (1987), pages 158 and 159], high breakdown point, high efficiency and a reliable algorithm to compute them. Moreover, it is desirable that their asymptotic theory satisfy two important features: be valid under verifiable regularity assumptions, and hold uniformly over a relatively large set of distribution functions with known size.

With these desired properties in mind we propose to use scale-equivariant $M$-location estimates calculated with a smooth function $\psi$ in the family (3) and with an $S$-scale estimate calculated with a function $\rho$ in Tukey's class (7). These $MM$-location estimates have simultaneously high breakdown point and high efficiency at the central model. Moreover, we showed that under verifiable conditions they are *uniformly consistent* and *uniformly asymptotically normal* over a contamination neighborhood of *known size*. For each choice of breakdown point and efficiency we showed how to compute the size of the contamination neighborhood where these uniform results hold. When the center of these neighborhoods is the standard normal distribution we found that these sizes range from 11% (for estimates with 50% breakdown point) to 25% (for 25% breakdown-point estimates). Hence, in most practical situations where the contamination is below 10% [Hampel, Ronchetti, Rousseeuw and Stahel (1986)], an $MM$-location estimate with 50% breakdown point has uniform asymptotic properties that allow for reliable statistical inference based on its asymptotic distribution.

It is of much interest to obtain this kind of uniform asymptotic properties for robust regression estimates. In principle, $MM$-regression estimates are good candidates to have satisfactory uniform asymptotic properties. Salibian-Barrera (2000) shows that under certain regularity conditions these estimates are asymptotically normal for any distribution in the contamination neighborhood. The main technical difficulty when using the approach presented in this paper to study $MM$-regression estimates seems to be to find sufficient regularity conditions on the loss function $\rho$ to show the uniform consistency of the $S$-regression estimate. Once this is established, Theorems 1 and 3 apply with appropriate modifications.

## APPENDIX

### A.1. Proofs.

PROOF OF THEOREM 1. For any $t \in \mathbb{R}$ and $F \in \mathcal{H}_\varepsilon$ let $\mu_\psi(t, F) = E_F \psi((X - t)/\boldsymbol{\sigma}(F))$, and fix an arbitrary $\tilde{\varepsilon} > 0$. Let $\boldsymbol{\sigma} = \boldsymbol{\sigma}(F)$, $\boldsymbol{\mu} = \boldsymbol{\mu}(F)$, $\psi(X, t, s) = \psi((X - t)/s)$. Also let $Y_i(t) = \psi(X_i, t, \hat{\sigma}_n)$ and $Y(F, t) = E_F \psi(X, t, \boldsymbol{\sigma})$. Let $\overline{\psi}_n(t) = \frac{1}{n} \sum_{i=1}^n Y_i(t)$ and $\mu_\psi(t, F) = E_F(\psi(X, t, \boldsymbol{\sigma}))$. For each $m \in \mathbb{N}$, $t \in \mathbb{R}$, $F \in \mathcal{H}_\varepsilon$ and $\tau > 0$ let

$$\mathcal{A}_m(F, t, \tau) = \left\{ \sup_{n \geq m} |\overline{\psi}_n(t) - \mu_\psi(t, F)| > \tau \right\};$$

then $\lim_{m\to\infty} \sup_{F\in\mathcal{H}_\varepsilon} P_F(\mathcal{A}_m(F,t,\tau)) = 0$. We have

$$\{\hat{\mu}_n < \boldsymbol{\mu} - \tilde{\varepsilon}\} \subseteq \{|\overline{\psi}_n(\boldsymbol{\mu} - \tilde{\varepsilon}/2) - \mu_\psi(\boldsymbol{\mu} - \tilde{\varepsilon}/2, F)| > a(\tilde{\varepsilon})\} = A_n(F, \tilde{\varepsilon}),$$

where $a(\tilde{\varepsilon})$ is given by $a(\tilde{\varepsilon}) = \inf_{F\in\mathcal{H}_\varepsilon} \mu_\psi(\boldsymbol{\mu}(F) - \tilde{\varepsilon}/2, F)$. Similarly

$$\{\hat{\mu}_n > \boldsymbol{\mu} + \tilde{\varepsilon}\} \subseteq \{|\overline{\psi}_n(\boldsymbol{\mu} - \tilde{\varepsilon}/2) - \mu_\psi(\boldsymbol{\mu} - \tilde{\varepsilon}/2, F)| > b(\tilde{\varepsilon})\} = B_n(F, \tilde{\varepsilon}),$$

where $b(\tilde{\varepsilon})$ equals $b(\tilde{\varepsilon}) = \inf_{F\in\mathcal{H}_\varepsilon} -\mu_\psi(\boldsymbol{\mu}(F) + \tilde{\varepsilon}/2, F)$. It is easy to see that $a(\tilde{\varepsilon}) = \inf_{F\in\mathcal{H}_\varepsilon} \mu_\psi(\boldsymbol{\mu}(F) - \tilde{\varepsilon}/2, F) > 0$, and that $b(\tilde{\varepsilon}) = \inf_{F\in\mathcal{H}_\varepsilon} -\mu_\psi(\boldsymbol{\mu}(F) + \tilde{\varepsilon}/2, F) > 0$. It follows that $\{|\hat{\mu}_n - \boldsymbol{\mu}| > \tilde{\varepsilon}\} \subseteq A_n(F, \tilde{\varepsilon}) \cup B_n(F, \tilde{\varepsilon})$. Then

$$\mathcal{M}_m(F, \tilde{\varepsilon}) = \left\{ \sup_{n\geq m} |\hat{\mu}_n - \boldsymbol{\mu}| > \tilde{\varepsilon} \right\}$$
$$\subseteq \mathcal{A}_m(F, \boldsymbol{\mu} - \tilde{\varepsilon}/2, a(\tilde{\varepsilon})) \cup \mathcal{A}_m(F, \boldsymbol{\mu} + \tilde{\varepsilon}/2, b(\tilde{\varepsilon})).$$

It follows that $\lim_{m\to\infty} \sup_{F\in\mathcal{H}_\varepsilon} P_F[\mathcal{M}_m(F, \tilde{\varepsilon})] = 0$. $\square$

PROOF OF THEOREM 2. We need to introduce the following notation. Let $\rho(x,t,s) = \rho((x-t)/s)$. Denote the set of positive real numbers $(0,\infty)$ by $\mathbb{R}_+$. For each $t \in \mathbb{R}$ and $s \in \mathbb{R}_+$ let $\gamma(F,t,s) = E_F \rho(X,t,s)$, and let $\gamma_n(t,s) = \gamma(F_n, t, s) = \frac{1}{n}\sum_{i=1}^n \rho(x_i, t, s)$, where $F_n$ denotes the empirical distribution function of the sample. It is easy to see that

$$\gamma(F, 0, \boldsymbol{\sigma}(F)) < \gamma(F, t, \boldsymbol{\sigma}(F)) \qquad \forall\, |t| \geq t^*.$$

Also, because of (16), there exists $\eta$ independent of $F$ such that

$$\inf_{-t^* \leq t \leq t^*, s^- \leq s \leq s^+} \gamma''(F,t,s) \geq \eta > 0 \qquad \forall\, F \in \mathcal{H}_\varepsilon,$$

where $\gamma''(F,t,s) = \partial^2 \gamma(F,t,s)/\partial t^2$. Hence the family of functions $\gamma(F, t, \boldsymbol{\sigma}(F))$ with $F \in \mathcal{H}_\varepsilon$ has a unique minimum in the fixed interval $(-t^*, t^*)$. For each $F \in \mathcal{H}_\varepsilon$ denote this unique minimum by $\tilde{\boldsymbol{\mu}}(F)$. Now fix an arbitrary neighborhood $B_\delta(\tilde{\boldsymbol{\mu}}(F))$ of $\tilde{\boldsymbol{\mu}}(F)$. Let $\tilde{\varepsilon}(\delta, F)$ satisfy

(19) $$\inf_{t \notin B_\delta(\tilde{\boldsymbol{\mu}}(F))} \gamma(F, t, \boldsymbol{\sigma}(F)) \geq \gamma(F, \tilde{\boldsymbol{\mu}}(F), \boldsymbol{\sigma}(F)) + \tilde{\varepsilon}(\delta, F).$$

We have that $\tilde{\varepsilon} = \tilde{\varepsilon}(\delta) = \inf_{F\in\mathcal{H}_\varepsilon} \tilde{\varepsilon}(\delta, F) > 0$. Choose an arbitrary $\tilde{\delta} > 0$ and let $I_2 \subset \mathbb{R}$ and $m_0 = m_0(\tilde{\delta})$ such that

(20) $$P_F[\tilde{\mu}_n \in I_2 \,\forall\, n \geq m] > 1 - \tilde{\delta} \qquad \forall\, m \geq m_0.$$

We can now build a finite coverage (independent of $F \in \mathcal{H}_\varepsilon$) of the set $I_2 \cap B_\delta(\tilde{\boldsymbol{\mu}}(F))^c$ with balls $B(t_1), \ldots, B(t_r)$ such that for every $j = 1, \ldots, r$ we have

$$E_F \left[ \inf_{t' \in B(t_j)} \rho(X, t', \boldsymbol{\sigma}(F)) \right] \geq \gamma(F, \tilde{\boldsymbol{\mu}}, \boldsymbol{\sigma}(F)) + \tilde{\varepsilon}.$$

For each of these centers $t_k$ let $Y_i(t_k) = \inf_{t' \in B(t_k)} \rho(X_i, t', \hat{\sigma}_n)$ and

$$Y(F, t_k) = E_F\left[\inf_{t' \in B(t_k)} \rho(X, t', \boldsymbol{\sigma}(F))\right] \neq E_F[Y_i(t_k)].$$

Consider the events

$$A_m(F, t_k) = \left\{\sup_{n \geq m} |\overline{Y}_n(t_k) - Y(F, t_k)| \leq \tilde{\varepsilon}\right\}, \qquad m \in \mathbb{N}.$$

There exists $m_1(\tilde{\delta})$ independent of $F$ such that

$$P_F(A_m(F, t_k)) > 1 - \tilde{\delta} \qquad \forall m \geq m_1(\tilde{\delta}), \ \forall F \in \mathcal{H}_\varepsilon, \ \forall t_k \in I_2.$$

Now note that

$$A_m(F, t_k) \subseteq \left\{\inf_{t \in B(t_k)} \frac{1}{n} \sum_{i=1}^n \rho(x_i, t, \hat{\sigma}_n) \geq \gamma(F, \tilde{\boldsymbol{\mu}}(F), \boldsymbol{\sigma}(F)) + 2\tilde{\varepsilon} \ \forall n \geq m\right\}$$
$$= C_m(F, t_k).$$

Let

$$D_m(F) = \left\{\frac{1}{n} \sum_{i=1}^n \rho(x_i, \tilde{\boldsymbol{\mu}}(F), \boldsymbol{\sigma}(F)) \leq \gamma(F, \tilde{\boldsymbol{\mu}}(F), \boldsymbol{\sigma}(F)) + \tilde{\varepsilon} \ \forall n \geq m\right\}.$$

There exists $m_2 = m_2(\tilde{\delta})$ (independent of $F$) such that for $m \geq m_2$ we have

$$P_F(D_m(F)) > 1 - \tilde{\delta} \qquad \forall F \in \mathcal{H}_\varepsilon.$$

Take $m_3 = \max(m_0, m_1, m_2)$. We have

$$P_F[C_m(F) \cap D_m(F)] \geq 1 - 2\tilde{\delta} \qquad \forall m \geq m_3, \ \forall F \in \mathcal{H}_\varepsilon.$$

We also have

$$C_m(F) \cap D_m(F) \subseteq [\tilde{\boldsymbol{\mu}}_m \in B_\delta(\tilde{\boldsymbol{\mu}}(F)) \ \forall m \geq m_2].$$

That is, for each $\delta > 0$ we have $\lim_{m \to \infty} \sup_{F \in \mathcal{H}_\varepsilon} P_F[\sup_{n \geq m} |\tilde{\boldsymbol{\mu}}_n - \tilde{\boldsymbol{\mu}}(F)| > \delta] = 0$. □

DEFINITION 2 (Uniform small $o$ in probability). Let $a_n$, $n \geq 1$, be a sequence of real numbers and let $X_n$, $n \geq 1$, be a sequence of random variables. We say that $X_n = Uo_P(a_n)$ over the set of distribution functions $\mathcal{H}_\varepsilon$ if, $\forall \delta > 0$,

$$\lim_{n \to \infty} \sup_{F \in \mathcal{H}_\varepsilon} P_F\left[\left|\frac{X_n}{a_n}\right| > \delta\right] = 0.$$

DEFINITION 3 (Uniformly asymptotically normal). We say that a sequence $X_n$, $n \in \mathbb{N}$, is uniformly asymptotically normal (UAN) over the set of distribution functions $\mathcal{H}_\varepsilon$ if

$$\sup_{F \in \mathcal{H}_\varepsilon} \sup_{x \in \mathbb{R}} |P_F(X_n \leq x) - \Phi(x)| = o(1). \tag{21}$$

With the above definitions we can show that these "uniform little $o$," "uniform big $O$" and "uniform asymptotic distribution" behave similarly to their "nonuniform" counterparts. In particular, if $a_n = Uo_P(1)$ and $X_n$ is UAN, then $X_n + a_n$ is UAN.

PROOF OF THEOREM 3. To simplify the notation, in what follows let $\mu = \mu(F)$, $\tilde{\mu} = \tilde{\mu}(F)$ and $\sigma = \sigma(F)$. The idea of the proof is to show that $\sqrt{n}(\hat{\mu}_n - \mu)$ can be represented as a linear term plus a uniformly small remainder. We use the Berry–Esseen theorem to show that the linear part is UAN (see Definition 3) and then argue that if $a_n = Uo_P(1)$ (see Definition 2) and $X_n$ is UAN, then $X_n + a_n$ is UAN. First note that by Theorem 2 and Remark 4 we have $\hat{\sigma}_n - \sigma = Uo_P(1)$, $\tilde{\mu}_n - \tilde{\mu} = Uo_P(1)$ and $\hat{\mu}_n - \mu = Uo_P(1)$. The idea of the proof is to show that

$$\sqrt{n} \frac{(\hat{\mu}_n - \mu)}{\sqrt{V}} = \sqrt{n} \frac{\overline{W}_n}{\sqrt{V}} + Uo_P(1), \tag{22}$$

where

$$W_i = \frac{(\psi((x_i - \mu)/\sigma) - d(\rho((x_i - \tilde{\mu})/\sigma) - b))}{e},$$

$$d = \frac{E_F\{\psi'((X - \mu)/\sigma)(X - \mu)/\sigma\}}{E_F\{\rho'((X - \tilde{\mu})/\sigma)(X - \tilde{\mu})/\sigma\}}, \tag{23}$$

$$e = E_F\left\{\psi'\left(\frac{(X - \mu)}{\sigma}\right)\right\}.$$

Note that $|W_i|$ are bounded and hence their moments are bounded uniformly for $F \in \mathcal{H}_\varepsilon$. The variance of $W_i$ is bounded away from zero uniformly on $F \in \mathcal{H}_\varepsilon$. The Berry–Esseen theorem yields

$$\sup_{F \in \mathcal{H}_\varepsilon} \sup_{x \in \mathbb{R}} \left| P_F\left\{ \frac{\sqrt{n}\,\overline{W}_n}{\sqrt{V}} < x \right\} - \Phi(x) \right| = o(1).$$

Now, to complete the proof, use that if $a_n = Uo_P(1)$ and $X_n$ is UAN, then $X_n + a_n$ is UAN. □

**Acknowledgment.** We thank two anonymous referees for their helpful comments that greatly improved the presentation of this paper.

UNIFORM ASYMPTOTICS FOR ROBUST LOCATION 15

School of Mathematics and Statistics  
Carleton University  
1125 Colonel By Drive  
Room 4302 Herzberg Building  
Ottawa, Ontario K1S 5B6  
Canada  
e-mail: matias@math.carleton.ca

Department of Statistics  
University of British Columbia  
Agricultural Road 333-6356  
Vancouver, British Columbia V6T 1Z2  
Canada  
e-mail: ruben@stat.ubc.edu